%% file: Main.tex
\documentclass[11pt,a4paper]{article}
\usepackage[top=1in, bottom=1.25in, left=1in, right=1in]{geometry}
\usepackage[utf8]{inputenc}
\usepackage[english]{babel}
\usepackage{graphicx,epstopdf}
\usepackage{caption, subcaption,float}
\usepackage{color}
\usepackage{multirow}
\usepackage{tikz}
\usepackage{amsthm}
\usepackage{tabls}
\usepackage{cite}
\usepackage[round,comma,sort]{natbib}

\usepackage[ruled, linesnumberedhidden,vlined]{algorithm2e} 
\usetikzlibrary{matrix}
\usepackage{amsmath,amssymb,dsfont,mathtools}
\usepackage{hyperref}
\hypersetup{
    colorlinks=true,       
    linkcolor=blue,          
    citecolor=olive,       
    filecolor=magenta,      
    urlcolor=blue          
}
\usepackage{aliascnt}

\usepackage{soul}

\begin{document}

\title{Pure infinitely braided Thompson groups}
\date{\today }
\author{Mar\'{i}a Cumplido}


\maketitle
\theoremstyle{plain}
\newtheorem{theorem}{Theorem}

\newaliascnt{lemma}{theorem}
\newtheorem{lemma}[lemma]{Lemma}
\aliascntresetthe{lemma}
\providecommand*{\lemmaautorefname}{Lemma}

\newaliascnt{proposition}{theorem}
\newtheorem{proposition}[proposition]{Proposition}
\aliascntresetthe{proposition}
\providecommand*{\propositionautorefname}{Proposition}

\newaliascnt{corollary}{theorem}
\newtheorem{corollary}[corollary]{Corollary}
\aliascntresetthe{corollary}
\providecommand*{\corollaryautorefname}{Corollary}

\newaliascnt{conjecture}{theorem}
\newtheorem{conjecture}[conjecture]{Conjecture}
\aliascntresetthe{conjecture}
\providecommand*{\conjectureautorefname}{Conjecture}

\newaliascnt{question}{theorem}
\newtheorem{question}[question]{Question}
\aliascntresetthe{question}
\providecommand*{\conjectureautorefname}{Question}

\theoremstyle{remark}

\newaliascnt{claim}{theorem}
\newaliascnt{remark}{theorem}

\newtheorem{claim}[claim]{Claim}
\newtheorem{remark}[remark]{Remark}
\newaliascnt{notation}{theorem}
\newtheorem{notation}[notation]{Notation}
\aliascntresetthe{notation}
\providecommand*{\notationautorefname}{Notation}

\aliascntresetthe{claim}
\providecommand*{\claimautorefname}{Claim}

\aliascntresetthe{remark}
\providecommand*{\remarkautorefname}{Remark}

\newtheorem*{claim*}{Claim}
\theoremstyle{definition}

\newaliascnt{definition}{theorem}
\newtheorem{definition}[definition]{Definition}
\aliascntresetthe{definition}
\providecommand*{\definitionautorefname}{Definition}

\newaliascnt{example}{theorem}
\newtheorem{example}[example]{Example}
\aliascntresetthe{example}
\providecommand*{\exampleautorefname}{Example}

\def\autorefspace{\hspace*{-0.5pt}}
\def\sectionautorefname{Section\autorefspace}
\def\subsectionautorefname{Section\autorefspace}
\def\subsubsectionautorefname{Section\autorefspace}
\def\figureautorefname{Figure\autorefspace}
\def\subfigureautorefname{Figure\autorefspace}
\def\tableautorefname{Table\autorefspace}
\def\equationautorefname{Equation\autorefspace}
\def\Itemautorefname{item\autorefspace}
\def\Hfootnoteautorefname{footnote\autorefspace}
\def\AMSautorefname{Equation\autorefspace}

\newcommand{\co}{\simeq_c}
\newcommand{\w}{\widetilde}
\newcommand{\po}{\preccurlyeq}
\newcommand{\dist}{\mathrm{d}}

\def\Z{\mathbb Z} 
\def\Ker{{\rm Ker}} \def\R{\mathbb R} \def\GL{{\rm GL}}
\def\HH{\mathcal H} \def\C{\mathbb C} \def\P{\mathbb P}
\def\SSS{\mathfrak S} \def\BB{\mathcal B} \def\PP{\mathcal P} 
\def\supp{{\rm supp}} \def\Id{{\rm Id}} \def\Im{{\rm Im}}
\def\MM{\mathcal M} \def\S{\mathbb S}
\newcommand{\bigveer}{\bigvee^\Lsh}
\newcommand{\wedger}{\wedge^\Lsh}
\newcommand{\veer}{\vee^\Lsh}
\def\diam{{\rm diam}}

\newcommand{\myref}[2]{\hyperref[#1]{#2~\ref*{#1}}}

\begin{abstract}
We show that pure subgroups of infinitely braided Thompson's are bi-orderable. For every finitely generated pure subgroup, we give explicit sets of generators. 

\medskip

{\footnotesize
\noindent \emph{2020 Mathematics Subject Classification.} 20F05; 20F36; 20F60; 20F65.

\noindent \emph{Key words.} Thompson groups; braid groups; pure braid groups; orderings; finitely generated groups.}

\end{abstract}


\section{Introduction}

Infinitely braided Thompson groups $BV_n(H)$ were introduced in \citep{ArocaCumplido2022} as a generalization of braided Thompson groups \citep{Dehornoy,Brin1,Brin2} and the Thompson groups $V_n$ \citep{Higman}. Since braids are torsion-free, braided Thompson groups and their generalized version are examples of Thompson-like groups without torsion, which is not the case for the traditional Thompson groups. The aim of this paper is to study the orderability and generators of their pure subgroups $BF_n(H)$, where $H$ is a subgroup of the pure braid group on $n$ strands (for detailed definitions see \autoref{section_def}). These groups, in the same way that $BV_n(H)$ are a combination of braid groups and Thompson groups $V_n$, can be seen as a mixture of pure braid groups and the Thompson groups $F_n$. Although it will not be explicit, the definitions of these groups arise from the concept of cloning system ---see \citep{Zaremsky, WitzelZaremsky, SZ2021, SZ2023}  for references about these constructions.--- 

In \citep{BurilloMeneses} the authors proved that all pure subgroups of braided Thompson groups admit a bi-ordering, by using the already pre-existing bi-ordering in the pure braid group. In \autoref{sectio_order}, we modify their arguments to prove that all $BF_n(H)$ are also bi-orderable (\autoref{theorem_ordering}). Since in bi-orderable groups the conjugates of positive elements are also positive, an immediate corollary is that $BF_n(H)$ cannot has generalized torsion, that is, there is no product of conjugates of an element that equals the identity. 

Recently, \cite{SkipperWu} proved that $BV_n(H)$ and $BF_n(H)$ are not only finitely presentable when $H$ is finitely presentable, but that it is of type ``$F_n$" \footnote{$F_n$ is the standard notation for this concept of finiteness and we use one time in this paper with quotation marks. Do not mistake with the Thompson group $F_n< V_n$.} if and only if $H$ is, which makes $BV_n(\mathcal{B}_n)$ and $BF_n(\mathcal{P}_n)$ of type ``$F_\infty$". However, for $BF_n(H)$, we only had explicit set of generators for the case $BF_2(\{1\})$ (\cite{Burillo}). In \autoref{section_gen}, we apply the techniques of \citep{ArocaCumplido2022} to give explicit sets of generators for every $BF_n(H)$ with a finitely generated $H$ (\autoref{Theorem_Gen2}). Our results lead us to a natural question: are the pure braided Thompson group $BF_n(\{1\})$ and the  pure infinitely braided Thompson group $BF_n(\mathcal{P}_n)$ is isomorphic? This is answered in the negative in  \autoref{conjecture_iso}.

\section{Previous Definitions}\label{section_def}

\subsection{Thompson groups}

We first recursively construct the $n$-adic Cantor set $\mathfrak{C}_n$ be the $n$-adic Cantor set: $\mathfrak{C}_n^1$ corresponds to first subdividing $\mathfrak{C}_n^0= [0,1]$ into $2n-1$ intervals of equal length, numbered $1, \ldots, 2n-1$ from left to right, and then taking the collection of odd-numbered subintervals. Then we renumber these intervals from left to right using labels $C^1_1,\dots, C^1_n$ and we obtain $\mathfrak{C}_n^2$ from $\mathfrak{C}_n^1$ by applying the same procedure to each interval $C^1_i$ to obtain $C^2_{(i-1)n +1}, C^2_{(i-1)n +2},\dots, C^2_{(i-1)n +n}$. In that way we can define $\mathfrak{C}_n^j$ for every $j$. $\mathfrak{C}_n$ is the intersection of all~$\mathfrak{C}_n^i$.

\smallskip

We can select a cover for~$\mathfrak{C}_n$ by choosing pairwise disjoint intervals of the form~$C_i^j$ from any~$\mathfrak{C}^j_n$. Now, given two covers $C$ and $C'$ with the same number of intervals, we define an affine and orientation preserving map from the elements of $C$ to the elements of~$C'$, and we restrict it to $\mathfrak{C}_n$. This restriction is a homeomorphism of~$\mathfrak{C}_n$. The Thompson's group~$V_n$ -- for a reference see \citep{Higman}-- is the set of all such homeomorphisms, and it is a group under composition. 

\smallskip

\begin{figure}[H]
\centering
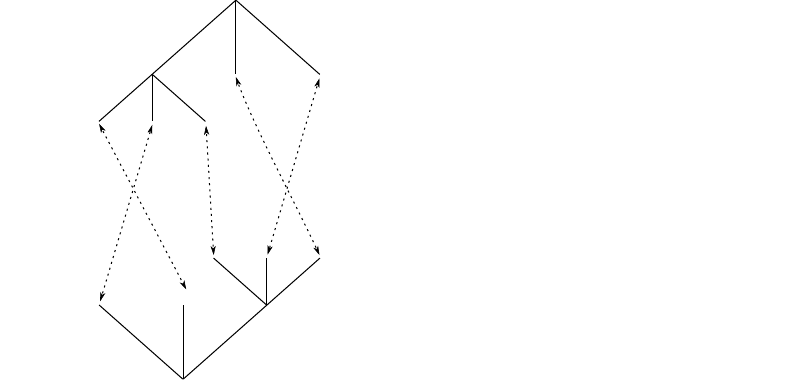
\caption{Two representatives of the same element in $V_3$.}
\label{thompsonequiv}
\end{figure}

It is easy to realise that the previous covers of $\mathfrak{C}_n$ are in bijection with finite full n-ary trees. An element of $V_n$ is then represented by a triple $(T, \tau, T')$, in which $T$ and $T'$ are full $n$-ary trees with the same number of leaves and $\tau$ is a bijection between those leaves. Notice that this representative is not unique (see \autoref{thompsonequiv}). We will work with this kind of representatives during this paper. A well known subgroup of $V_n$, that will be used in the last section, is $F_n$, consisting of all elements represented by triples $(T, \tau, T')$ where $\tau$ is trivial, so one can just write~$(T,T')$.

Since representatives in these groups are not unique, we need to stablish some specific terminology to define equivalence relations. A subtree of an $n$-ary $T$ formed by $n$ leaves sharing a parent plus this parent and the edges from the leaves to the parent is called a \emph{caret}. Starting by 1, enumerate from left to right the leaves of $T$. Attaching to the $i$-th leaf of $T$ a caret is denoted $T[i]$ (see \autoref{thompsonequiv}). 

\subsection{Braid groups and pure braid groups}

\begin{figure}[H]
     \begin{subfigure}[b]{0.5\textwidth}
         \centering
         \includegraphics[width=0.5\textwidth]{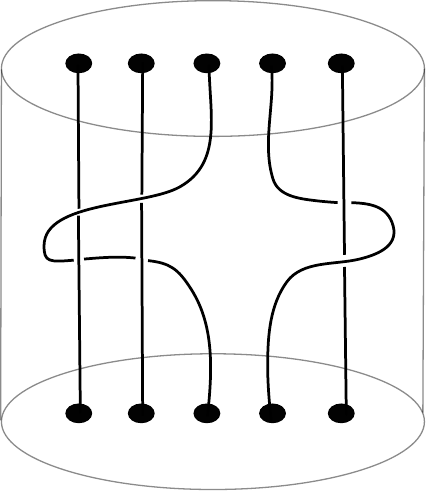}
         \caption{Braid as an object in a cylinder.}
         \label{trenzacilindro}
     \end{subfigure} \hfill
          \begin{subfigure}[b]{0.5\textwidth}
         \centering
         \includegraphics[width=0.3\textwidth]{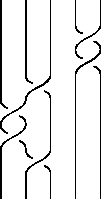}
         \caption{A standard projection of the braid.}
         \label{trenzaproy}
     \end{subfigure}
\caption{The pure braid $\sigma_4^2 \sigma_2\sigma_1^2\sigma_2$.}\label{trenza}
\end{figure}

A braid on $n$-strand is an  object composed of a cylinder with $n$ marked points in both its superior and inferior disk, that are linked by $n$-strands running monotonically in the vertical direction that never touch each other (see \autoref{trenzacilindro}). We say that two braids are equivalent if one can transform one into the other by continuous deformation without crossing strands. Usually, braids are graphically represented using a projection on the plane using crossings of two consecutive strands (see \autoref{trenzaproy}). To compose two braids, one glues the two corresponding cylinders and re-escalate, which is equivalent to concatenate both corresponding diagrams in the obvious way. The class of equivalence classes of braids with this composition is the braid group on $n$ strands, $\mathcal{B}_n$, introduced by \citep{Artin} and its classic presentation is the following one:
$$\mathcal{B}_n=\langle \sigma_1,\dots,\sigma_{n-1} \, |\, \sigma_i\sigma_j=\sigma_j\sigma_i  \text{ if } |i-j|>1; \sigma_i\sigma_j\sigma_i=\sigma_j\sigma_i\sigma_j  \text{ if } |i-j|=1 \rangle, $$
where $\sigma_i$ is the braid in which the strand in the position $i$ passes under the strand in the position $i+1$, and $\sigma_i^{-1}$ is the braid in which the strand in the position $i$ passes over the strand in the position $i+1$.

The pure braid group on $n$ strands, $\mathcal{P}_n$,  is the subgroup of $\mathcal{B}_n$ made of braids in which every strand starts and ends in the same position, that is, the $i$-th marked point in the upper disk of the cylinder is linked with the $i$-th marked point in the lower disk. The braid in \autoref{trenza} is a pure braid.

\begin{figure}[H]
\centering
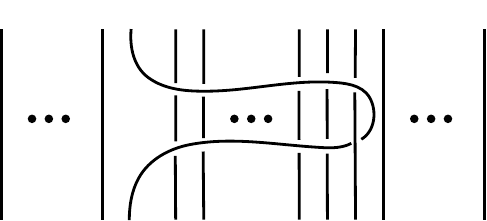
\caption{The pure generator $A_{i,j}$.}
\label{pure_gen}
\end{figure}

The pure braid group on $n$ strands has also a finite presentation \citep[Lemma~1.8.2]{Birman}, but in this paper we are going to use the set of finite generators (see \autoref{pure_gen}), which are all braids that can be written as $$A_{i,j}:= \sigma_i^{-1} \cdots \sigma_{j-2}^{-1}\sigma_{j-1}^{-2}\sigma_{j-2}\cdots \sigma_i, \quad 1\leq i <j \leq n. $$

\begin{figure}
\centering
\includegraphics[width=0.3\textwidth]{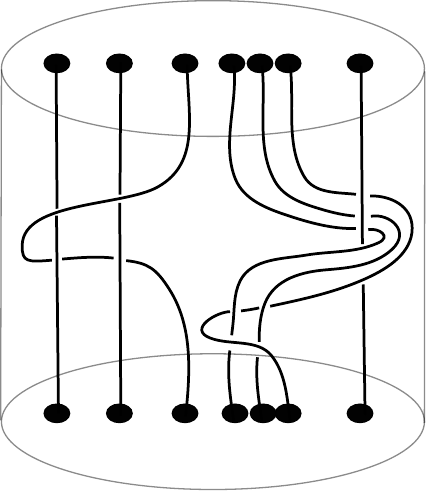}
\caption{$b[4,\sigma_2^{-1}\sigma_1^{-2}\sigma_2^{-1}]$, where $b$ is the braid from \autoref{trenza}.}
\label{trenzacilindro2}
\end{figure}

The groups that we will study mix Thompson groups and braid groups. To define equivalence relations in those groups, we will need the notion of \emph{splitting and re-braiding}. Fix some $n\in \mathbb{N}$. Given two braids $b$ (with any number of strands) and $\ell\in \mathcal{B}_n$ we denote $b[i,\ell]$ the braid obtained from $b$ by splitting in $n$ strands the strand in the $i$-th position and braiding these $n$ strands as~$\ell$ indicates (see an example in \autoref{trenzacilindro2}). 

\subsection{Infinitely braided Thompson groups $BV_n(H)$}

Infinitely braided Thompson groups\citep{ArocaCumplido2022} are, as presented in their name, a generalization of the well known families of braided Thompson groups introduced separately by \cite{Brin1,Brin2} and \cite{Dehornoy}.
Here we give new definitions that are more compact that the ones in the previous papers. The reader can easily check that they are in fact equivalent to the definitions in \citep[Section~2]{ArocaCumplido2022}.

\begin{figure}[h]
     \begin{subfigure}[a]{0.5\textwidth}
         \centering
         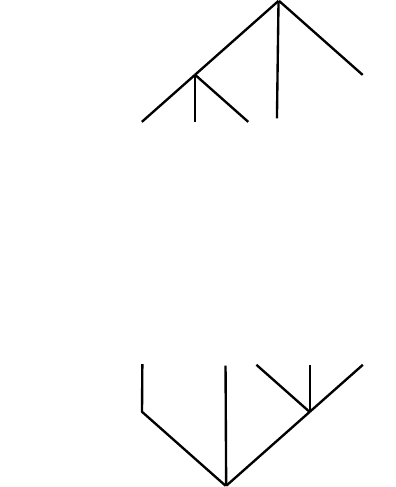
         \caption{The element $x=(T_1,\sigma_4^2 \sigma_2\sigma_1^2\sigma_2,\lambda,T_2)$, where \\ $\lambda_4=\sigma_2^{-1}\sigma_1^{-2}\sigma_2^{-1}$ and $\lambda_i$ is trivial for $i\neq 4$.}
         \label{generalizedelement}
    \end{subfigure}
       \begin{subfigure}[a]{0.5\textwidth}
         \centering
         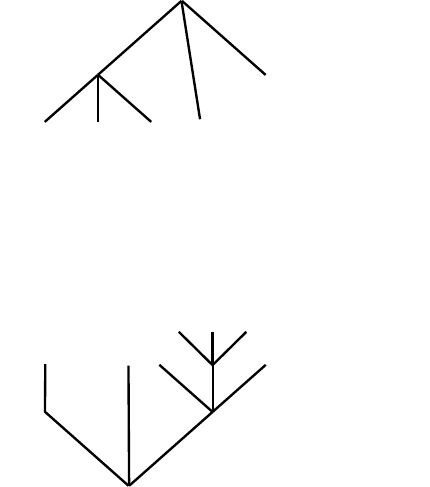
         \caption{The element $\alpha'_{3,4}(x)$ (trivial labels are omitted).}
         \label{generalizedelement2}
    \end{subfigure}
\caption{Examples of elements in $BV_3(\mathcal{B}_3)$.}\label{genel}
\end{figure}

The elements of $BV_n(H)$ with $H<\BB_n$ are obtained from sets $C_{m,n}$, of elements having the form $x=(T_1,b,\lambda,T_2)$, where $T_1$ and $T_2$ are $n$-ary trees with $m$ leaves, $b\in \BB_m$, and $\lambda$ is a list $\{\ell_,\dots,\ell_{m}\}$, $\ell_i\in H,\,  1\leq i < m$. Graphically we can see these elements as two trees whose leaves are connected by the braid $b$ that has the label $\ell_i$ in its $i$-th strand (see \autoref{generalizedelement}). Notice that, if $\lambda_1=\{\ell_1,\dots,\ell_{m}\}$ and $\lambda_2=\{\ell'_1,\dots,\ell'_{m}\}$, we can compute $\lambda_1\lambda_2=\{\ell_1\ell_1',\dots,\ell_{m}\ell_{m}'\} $. Moreover, we define $b(\lambda)=\{\ell_{b(1)},\dots,\ell_{b(m)}\}$, where $b(i)$ is the position of the point in the lower disk of $b$ link with the point $i$ in its upper disk. 

Associated with these sets, we define maps $\alpha'_{m,i} : C_{m,n} \rightarrow C_{m+n-1,n}$, $1\leq i < m$, such that $$\alpha'_{m,i} (T_1,b,\{\ell_1,\dots,\ell_i,\dots,\ell_m\},T_2)=(T_1[i],b[i,\ell_i],\{\ell_1,\dots,\underbrace{\ell_i,\dots,\ell_i}_n,\dots,\ell_{m+n-1}\},T_2[i]).$$ Given an element $x$, we say that $y$ is an \emph{expansion} of $x$ if it is obtained from $x$ by several $\alpha'$-applications or, alternatively, we say that that $x$ is a  \emph{reduction} of $y$. We also say that we obtain $y$ by \emph{expanding} $x$ or that we obtain $y$ by \emph{reducing} $x$.

\begin{figure}[h]
\centering
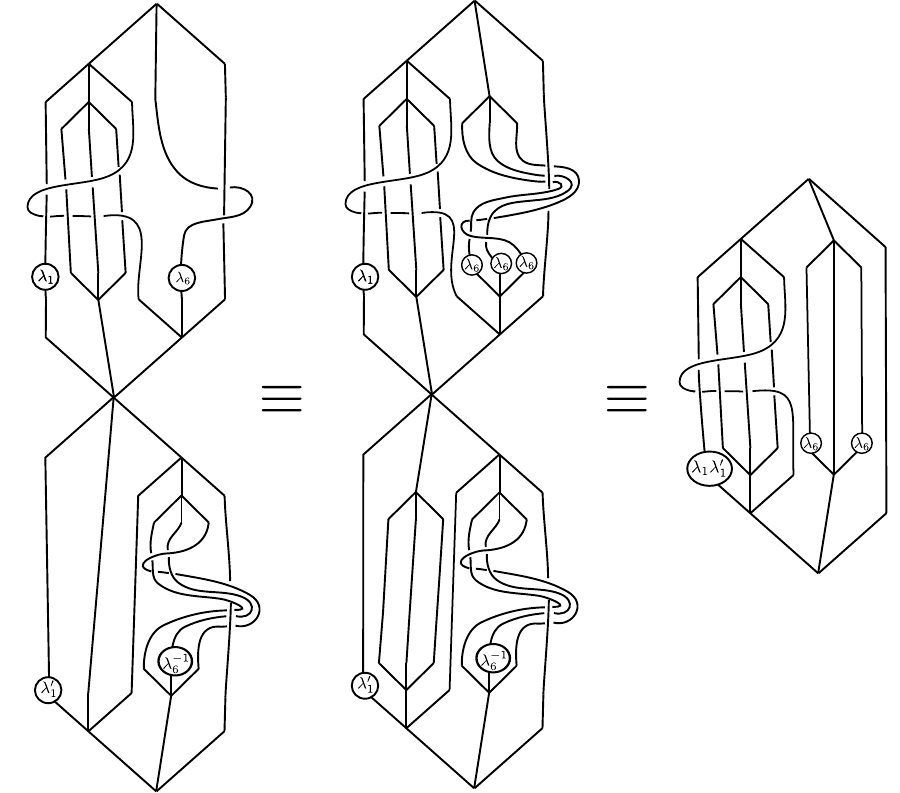
\caption{Composing two elements in $BV_3(\mathcal{B}_3)$. Trivial labels are omitted and $\lambda_6= \sigma_2^{-1}\sigma_1^{-2}\sigma_2^{-1}$.}
\label{composition}
\end{figure}

To compose two elements $x=(T_1,b_1,\lambda_1,T_2)$ and $y=(T_3,b_2,\lambda_2,T_4)$, one selects an expansion $(T_1',b_1',\lambda'_1,T_2')$ of $x$ and an expansion $(T_3',b_2',\lambda_2',T_4')$ of $y$ so that $T_2'=T_3'$ and $xy:=(T_1', b_1'b_2',b_2'(\lambda_1)\lambda_2, T_4' )$. The reader can see a graphic example in \autoref{composition}.

The infinitely braided Thompson group is given by the previously explained composition and the following infinite presentation: 
$$BV_n(H)= \langle x\in C_{m,n} \,|\,    x=\alpha'_{m,i} (x), \,  1 \leq i \leq n,\, m>0 \rangle.$$

\begin{remark}
Following the references given in this section, we remark that the same definitions work when taking forests instead trees by doing the obvious generalizations (as in the construction of Higman-Thompson's groups described in \citep{Higman}). For simplicity, in this paper we will deal only with trees, but similar techniques work for forests. 
\end{remark}

\subsection{Pure infinitely braided Thompson groups $BF_n(H)$}

In this article we will focus in a specific kind of subgroup of infinitely braided Thompson's groups, that this time mixes the pure braid group and the group $F_n$. This a generalization of the pure braided Thompson group defined in \citep{Burillo}.
The elements of $BF_n(H)$ with $H<\PP_n$ will be obtained from sets~$P_{m,n}$, of elements of the form $x=(T_1,p,\lambda,T_2)$, where $T_1$ and $T_2$ are $n$-ary trees with $m$ leaves, $b\in \PP_m$, and $\lambda=\{\ell_1,\dots,\ell_{m}\}$, $\ell_i\in H,\,  1\leq i < m$. Clearly, every~$P_{m,n}$ is a subset of~$C_{m,n}$. Also~$\alpha'_{m,i}$ restricts to~$\alpha_{m,i}: P_{m,n}\rightarrow P_{m+n-1,n}$, because~$p[i,\ell_i]$ is still a pure braid.

$$BF_n(H)= \langle x\in P_{m,n} \,|\,    x=\alpha_{m,i} (x), \,  0 \leq i \leq n,\, m>0 \rangle.$$

Since $\PP_n<\BB_n$, it is easy to check that~$BF_n(H)$ satisfies the properties to be a subgroup of~$BV_n(H)$. At this point the reader can realise that the elements pictured in \autoref{genel} and \autoref{composition} are elements in~$BF_3(\mathcal{P}_3)$.

\section{Bi-orderability of $BF_n(H)$}\label{sectio_order}

The first aim of this paper is to prove that the pure infinitely braided Thompson group $BF_n(H)$ is always bi-orderable for any $H$. We will do that by adapting the arguments of \citep{BurilloMeneses}. In their article, they construct a bi-ordering on the pure braided Thompson group, which is isomorphic to $BF_2(\{1\})$.  This ordering is based on the bi-ordering defined for pure braids, that we briefly describe next (all details can be found in \citep{KimRolfsen}).

\medskip

There is a bi-ordering on the free group $\mathcal{F}_n$ called the Magnus ordering. Consider the ring $\mathbb{Z}[X_1,\dots,X_n]$ of $n$ non-commutative variables with coefficients in $\mathbb{Z}$. There is a total order on the set of monomials of $\mathbb{Z}[X_1,\dots,X_n]$, since we can order the variables $X_1<\dots<X_n$ and say that a monomial is bigger than a monomial of less degree and use the lexicographical order if the degrees coincide. Now, if $\mathcal{F}_n$ is freely generated by $x_1,\dots, x_n$, we can define a homomorphism $\phi:\mathcal{F}_n\rightarrow \mathbb{Z}[X_1,\dots,X_n]$ such that $\phi(x_i)= 1+ X_i$. Define an element $x\in \mathcal{F}_n$ as positive if the coefficient of the smallest monomial of $\phi(x)$ is positive with respect the standard order in~$\mathbb{Z}$. This defines a bi-order for $\mathcal{F}_n$. The previous ordering induces a bi-order of $\mathcal{P}_n$ by using the lexicographical order in the well-known decomposition $$\mathcal{P}_n=(\cdots ((\mathcal{F}_1 \ltimes \mathcal{F}_2 )\ltimes \mathcal{F}_3)\dots) \ltimes \mathcal{F}_{n-1}),$$
where the semi-direct product comes from the splitting short exact sequence
$$ 1 \longrightarrow \mathcal{F}_n \longrightarrow \mathcal{P}_n \longrightarrow \mathcal{P}_{n-1} \longrightarrow 1.$$ Here $\PP_n$ surjects into $\mathcal{P}_{n-1}$ by just erasing its first strand. The sequence obviously splits by adding one trivial strand to the left of the braids in $\mathcal{P}_{n-1}$. Since the action of $\mathcal{P}_{n-1}$ in $\mathcal{F}_n$ preserves the order in $\mathcal{F}_n$, we can use the lexicographical order to create an order for the whole decomposition  \citep[Lemma~1]{KimRolfsen}.

\medskip

We will also need to have a bi-ordering on the Thompson group $F_n$. We can easily define a bi-order for $F_n$ based on the bi-order defined for $F_2$ in \citep{Dehornoy}: recall that an element~$f$ of $F_n$ is an order-preserving bijection between subintervals of $[0,1]$, each having length $\frac{1}{n^k}$. Take the first pair of corresponding intervals of different length and divide the length of the first by the length of the second (the slope). We say that  $f$ is positive if this number is a positive power of $n$. 

\begin{proposition}\label{Prop_JuanPep}
The group $BF_n(\{1\})$ is bi-orderable.
\end{proposition}

\begin{proof}
The proof works exactly the same as for $n=2$ \citep{BurilloMeneses}. We explain the scheme here for completeness and a better understanding of the proofs for the general case. First, we consider the subgroup $PVB$ of $BF_n(\{1\})$ which elements can be written as $(T,p,\lambda_{id},T)$, where $\lambda_{id}$ is a list formed only by copies of $1$. That is, we consider the elements of~$BF_n(\{1\})$ in which both trees are the same ---notice that $PVB$ is a subgroup because this is consistent with the equivalence relation of the group---. We say that $(T,p,\lambda_{id},T)$ is positive if~$p$ is positive with respect to the order described above. In \citep[Lemma~3.1]{BurilloMeneses}, the authors prove that this definition of positivity is consistent with the equivalent relation of the group for $n=2$, that is, if we split a strand of a positive braid~$p$, then the resulting braid is still positive. This argument works for any~$n$, as splitting a strand in~$n$ new strands is just splitting in two strands $n$~times. After this, it is proven in \citep[Corollary~3.2]{BurilloMeneses} that the set of positive elements~$\mathcal{P}$ in~$PVB$ defines a bi-order, by showing that $\mathcal{P}$ is a semigroup invariant under conjugation and that $PVB=\mathcal{P}\sqcup \{1\} \sqcup \mathcal{P}^{-1}$. This proof works for any $n$.  

Finally, \citep[Theorem~3.3]{BurilloMeneses} shows that this induces an order in $BF_2(\{1\})$, and this generalizes without any changes for any~$n$: observe that there is a natural splitting short exact sequence 

$$ 1 \longrightarrow PVB \longrightarrow BF_n(\{1\}) \longrightarrow F_n \longrightarrow 1,$$
where the maps from $F_n$ to $BF_n(\{1\})$ just send $(T_1,T_2)$ to $(T_1,1,\lambda_{id}, T_2)$. This gives us the decomposition $BF_n(\{1\})= PVB \ltimes F_n$.  To see that the bi-order in $PVB$ and the bi-order in $F_n$ induce a bi-order in $BF_n(\{1\})$ one just need to show that the action of $F_n$ in $PVB$ preserves the bi-order on $PVB$ \citep[Lemma~1]{KimRolfsen}, which is proven by realizing that the action of an element of $F_n$ on an element $(T,p,\lambda_{id}, T)$ is just an equivalent element $(T',p',\lambda_{id}, T')$ where $p'$ is obtained from $p$ by splitting or merging strands, so it maintains the positivity. 

\end{proof}

We want to prove that $BF_n(H)$ is bi-orderable for any $H$. Our proof will mimic the previously explained strategy of \citep{BurilloMeneses} with additional steps to deal with labels.
First we define as $PVB(H)$ the set of elements $(T,p, \lambda, T)\in BF_n(H)$ under the same equivalence relation as in $BF_n(H)$, which is indeed a group.

\begin{theorem}\label{theorem_PVB}
The group $PVB(H)$ is bi-orderable. 
\end{theorem}

\begin{proof}
 
We will say that an element is positive $(T,p,\{\ell_1,\dots,\ell_m\}, T)$ if the first non-trivial $\ell_i$ is positive, or if every $\ell_i$ is trivial and $p$ is positive with the usual bi-order for pure braids. That is, we use the lexicographic order. 

First we need to prove that this definition is consistent with the equivalence relation of the group. For this, just notice that the splitting of strands does not change the positivity of the first non-trivial label. Also, if every label is trivial, then this has been previously proven in \autoref{Prop_JuanPep}.

Now, we will need to show that the set $P$ of positive elements is a semigroup such that $PVB(H)=P\sqcup \{1\} \sqcup P^{-1}$.  To check that it is a semigroup, take two positive elements $(T,p_1,\{\ell_1,\dots,\ell_m\},T)$ and $(T,p_2,\{\ell'_1,\dots,\ell'_m\},T)$ ---notice that we can suppose that the have the same tree $T$ because we have proven that splitting of strands does not change positivity, so we can expand the elements as we want---. Their product is  $(T,p_1p_2,\{\ell_1\ell_1',\dots,\ell_m\ell_m'\},T)$ and, since the first non-trivial $\ell_i$ and $\ell'_i$ are positive, then the first non trivial $\ell_i\ell'_i$ is also positive. The case where all labels are trivial has been already proven in \autoref{Prop_JuanPep}. Therefore, $(T,p_1p_2,\{\ell_1\ell_1',\dots,\ell_m\ell_m'\},T)\in P$  and $P$ is a semigroup. The inverse of $(T,p_1,\{\ell_1,\dots,\ell_m\},T)$ is $(T,p_1^{-1},\{\ell_1^{-1},\dots,\ell_m^{-1}\},T)$ so it follows immediately that $PVB(H)=P\sqcup \{1\} \sqcup P^{-1}$. Finally, we check that $P$ is invariant under conjugation. Conjugate a positive $(T,p_1,\{\ell_1,\dots,\ell_m\},T)$ by $(T,p_2,\{\ell'_1,\dots,\ell'_m\},T)$ ---again, we can suppose that the trees are the same---. The result of this conjugation is $$(T,p_2^{-1}p_1p_2,\{\ell_1'^{-1}\ell_1\ell'_1,\dots,\ell_m'^{-1}\ell_m\ell_m'\},T).$$
 We know that the set of positive pure braids is invariant under conjugation, so positivity is maintained. 

\end{proof}

\begin{theorem}\label{theorem_ordering}
The group $BF_n(H)$ is bi-orderable. 
\end{theorem}

\begin{proof}
We have the following short exact sequence: 
$$ 1 \longrightarrow PVB(H) \longrightarrow BF_n(H) \longrightarrow F_n \longrightarrow 1,$$
where $PVB(H)$ naturally injects in $BF_n(H)$ and $(T_1,p,\lambda, T_2)\in BF_n(H)$ is sent to $(T_1,T_2)\in F_n$. This sequence splits by sending $(T_1,T_2)$ to $(T_1,1,\lambda_{id}, T_2)$.  Now we have to prove that the action of $F_n$ in $PVB(H)$ preserves the set $P$ of positive elements in $PVB(H)$. If that is true, then by \citep[Lemma~1]{KimRolfsen}, we can use the lexicographic order on $BF_n(H)=PVB(H) \ltimes F_n$ induced by the previously described bi-ordering on $PVB(H)$ and the bi-ordering on $F_n$. The action by conjugation of an element $(T_1,T_2)\in F_n$ on an element $(T,p,\lambda,T)\in PVB(H)$ is $(T',p',\lambda',T')$, where $p'$ is obtained from $p$ by strand splitting and $\lambda'$ is obtained from $\lambda$ by replacing labels $\ell_i$ by a set of copies of it $\ell_i,\dots,\ell_i$. Therefore, the first non-trivial label is still positive or, if all labels are trivial, $p'$ is positive (\autoref{Prop_JuanPep}). 
\end{proof}

\begin{remark}\label{remark_Ishida}
It is worth mentioning the work of \cite{Ishida} that studies bi-orderings for cloning systems. In our particular context, we have groups $P_m \times P_n^m$ and every expansion defines a cloning operation from $P_m \times P_n^m$ to $P_{m+n-1} \times P_n^{m+n-1}$. The author proves that $BF_n(H)$ is bi-orderable if we can put a bi-order on $P_m \times P_n^m$ that is maintained under the cloning operation, which is exactly what we prove in \autoref{theorem_PVB}. The proofs in \citep{Ishida} are made for $n=2$ but they extend to $n>2$ using the concept of $n$-cloning system \citep{SZ2023}.  
\end{remark}

\section{Explicit generating set for $BF_n(H)$}\label{section_gen}

In this section we want to provide an explicit finite set of generators for $BF_n(H)$ when $H$ is finitely generated. More specifically, we will give an upper bound for the minimal number of needed generators and we will describe them. For $BF_2(\{1\})$, \cite{Burillo} were able to give explicit infinite presentations and a finite presentation with 10 generators (and 192 relations) based on the infinite and finite presentations of $BV_2\{1\}$ in the thesis of \cite{belkthesis} and the finite presentation for the pure braid group \citep[Lemma~1.8.2]{Birman}. Their approach describes in detail how every generator of an infinite presentation decomposes as a product involving powers of the 10 purposed generators. In this article, we are going to replicate the techniques used in \citep{ArocaCumplido2022} to give an explicit set of generators for every $BV_n(H)$ with a finitely generated $H$, based on the arguments of \cite{Higman} to prove the finite generation of $BV_n(\{1\})$. The idea is to decompose every element as a product of elements whose trees have as less leaves as possible, and those elements will be our generators. 

\begin{definition}
We say that a generator $A_{i,j}\in \mathcal{P}_m$ is \emph{$n$-irreducible} if $i\leq n$, $j-i \leq n$ and $m-j<n$. Otherwise, we say that $A_{i,j}$ is \emph{$n$-reducible}. Intuitively, a $n$-irreducible $A_{i,j}$ cannot have $n$ consecutive strands that could be glued to one single strand when being part of an element in $BF_n(\{1\})$ and performing a reduction.

\end{definition}

\begin{proposition}\label{Theorem_Gen1}
The group $BF_n(\{1\})$ can be generated by $n^2+2n+1$ generators for $n\neq 2$ and 10 generators for $n=2$.  In particular, it can be generated by the following elements:

\begin{itemize}

\item The $n$ generators equivalent to the ones of $F_n$ described by \cite{Brown}, that is, the elements $$\left(R[n],1,\lambda_{id}, R[i]\right),$$ for $i=1,\dots, n-1$ and $$\left(R[n][2n-1], 1,\lambda_{id}, R[n][n]\right),$$ where $R$ is the full $n$-ary tree with one caret. 

\item $8$ (for $n=2$) or $n^2+n+1$ (for $n> 2$) generators having the form $(T,A_{i,j},\lambda_{id},T)$, where~$A_{i,j}$ is a $n$-irreducible generator (in some $\mathcal{P}_m$) and $T$ can be any full $n$-ary tree with $m$ leaves. One can choose a different $T$ for each generator.

\end{itemize}

\end{proposition}

\begin{figure}[h]
\centering
\includegraphics[scale=0.8]{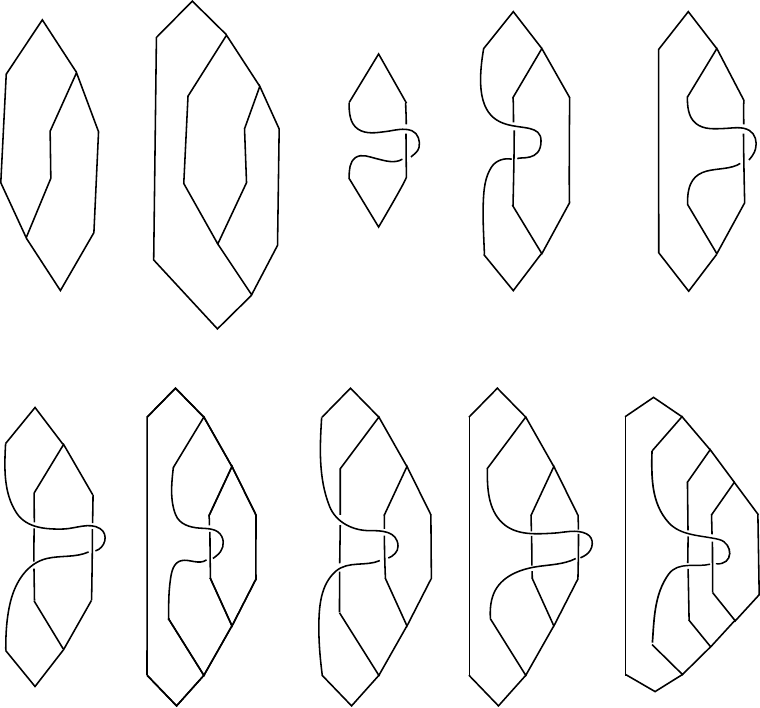}
\caption{Generators of $BF_2(\{1\})$.}\label{Fig_gen1}
\end{figure}

\begin{remark}
For $n=2$, if one takes in the previous result $T$ as the tree obtained from $R$ by attaching recursively carets to the last leave the suitable number of times, we obtain exactly the 10 generators given in \citep[Theorem~6.1]{Burillo}, which are depicted in \autoref{Fig_gen1}. 
\end{remark}

\noindent 
\begin{proof}[Proof of \autoref{Theorem_Gen1}]
First notice that, given any element $x=(T_1,p,\lambda_{id},T_2)$, we can choose any $T$ with the same leaves as $T_1$ (or $T_2$) and write $$x=(T_1,1,\lambda_{id},T)(T,p,\lambda_{id},T)(T,1,\lambda_{id},T_2).$$ Any element $(T',1,\lambda_{id},T'')$ can be written using the $n$ equivalent generators to the ones of~$F_n$ given by \cite{Brown}. Thus, the problem reduces to find a finite set of elements that allow us to generate the set of element of the form $(T,p,\lambda_{id},T)$.

Now, we know that $p$ can be written using the generators of pure braid groups (see \autoref{section_def}), so $(T,p,\lambda_{id},T)$ is a product of elements having the form $(T,A_{i,j},\lambda_{id},T)$. By the previous paragraph, for every $(T,A_{i,j},\lambda_{id},T)$ we can choose $T$ to be any $n$-ary tree where its number of leaves~$m$ is such that $A_{i,j}\in \mathcal{P}_m$. In particular,  if $A_{i,j}$ is $n$-reducible, one could choose a $T$ such that $$(T,A_{i,j},\lambda_{id},T)\equiv \alpha_{m-n+1,k}(T',A_{i',j'},\lambda_{id},T'),$$ for some $k$. That is, there are $n$ consecutive strands in $A_{i,j}$ that will be glued to one single strand when applying a reduction. Iterating this process, we arrive to the conclusion that, besides the $n$ generators previously mentioned, we just need the elements $(T,A_{i,j},\lambda_{id},T)$, where $A_{i,j}$ is $n$-irreducible to generate the whole group. To finish the proof, we need to study  the cardinality of this set of generators.

Notice that $m$ cannot be any number. In fact $m\in\{1,n,2n-1,3n-2,4n-3, 5n-3\dots\}$. It is easy to check that if $n=2$, any $A_{i,j}\in \mathcal{P}_m$ will be $n$-reducible for $m>4n-3$, and, if $n>2$, the same stands for $m\geq 4n-3$. We are going to count how many irreducible pure generators there are  for each $m$, by taking into account that we have three sets of strands whose number cannot exceed $n-1$: the set $S_1$ of strands before the position $i$, the set $S_2$ between the position $i$ and $j$ (without taking into account $i$ and $j$), and the set $S_3$ after the position $j$. For $m=n$, we have $n-2$ strands to distribute between the sets $S_1, S_2$ and $S_3$, which reduces to see which are the numbers for $|S_1|$ and $|S_2|$ satisfying $|S_1|+|S_2|\leq n-2$. This gives us $\dfrac{(n-1)n}{2}$ possibilities. For $m=2n-1$, we have $2n-3$ strands to distribute between the sets $S_1, S_2$ and~$S_3$. Note that in this case we have the restrictions $n-2\leq |S_1|+|S_2| \leq 2n-3$. The first inequality tells us that if $|S_1|=0$, then we have two possibilities for $|S_2|$; if $|S_1|=1$, then we have three possibilities for $|S_2|$, and so on till $|S_1|\in\{n-1,n-2\}$ that we have $n$ possibilities for $|S_2|$. Also, since $(n-1)+(n-1)=2n-2$, the second inequality only disallows one combination. Then we have $\dfrac{(n+1)(n+2)}{2}-3$ possibilities in total for this case. For $m=3n-2$, we have $3n-4$ strands to distribute between the sets $S_1, S_2$ and $S_3$. Note that in this case imposes the restrictions $2n-3\leq |S_1|+|S_2| \leq 3n-4$, which gives us just 3 possibilities. Adding the three numbers we obtain $n^2+n+1$. Finally, in the case $n=2,m=5$, we have just one possibility. 
\end{proof}

\begin{theorem}\label{Theorem_Gen2}
The group $BF_n(H)$, where $H$ is generated by $k$ generators, can be generated by $n^2+(k+2)n+1$ generators if $n\neq 2$ and $10+2k$ generators if $n=2$. In particular, it can be generated by the generators of $BF_n(\{1\})$ given in \autoref{Theorem_Gen1} and the $nk$ generators that can be written as $(R,1,\lambda_{i,h},R)$, for $1\leq i\leq n$, where $R$ is the full $n$-ary tree with one caret, and $\lambda_{i,h}=\{\ell_1,\dots,\ell_n\}$ is such that $\ell_i=h$ is a generator of $H$ and $\ell_j=1$ if $j\neq i$.
\end{theorem}

\begin{proof}

First notice that $(T_1,p,\{\ell_1,\dots,\ell_n\}, T_2)\equiv(T_1,p,\lambda_{id},T_2)(T_2,1,\{\ell_1,\dots,\ell_n\},T_2)$. By \autoref{Theorem_Gen1}, we already know how to generate the first factor. Also, the second factor is equivalent to a finite product of elements having the form $(T_2,1,\lambda_{i,h},T_2)$, for $0<i\leq m$, where $m$ is the number of leaves of $T_2$, and for every $T$ with the same number of leaves as $T_2$ we have that $$(T_2,1,\lambda_{id},T)(T,1,\lambda_{i,h},T)(T,1,\lambda_{id},T_2).$$ Therefore, we are going to prove that we can generate $(T,1,\lambda_{i,h},T)$, for a certain $T$ having any (possible) number of leaves, using the generators of the statement. 

If we want $i$ to be bigger than $n$, by expanding the first strand $q$ times, we can obtain from $(R,1,\lambda_{i,h},R)$, with $1<i\leq n$, a $(T,1,\lambda_{i+q(n-1),h},T)$ where $T$ has $(q+1)n-q$ leaves. Otherwise, expand the last leaf. This finishes the proof.
\end{proof}

In our last theorem, we use our generators containing non-trivial labels to reduce the set of generators for $BF_n(\mathcal{P}_n)$.

\begin{theorem}\label{Theorem_Gen3}
The group $BF_n(\mathcal{P}_n)$ can be generated by  9 generators for $n=2$ and $\frac{n^3}{2}+\frac{3n}{2}+1$ generators for $n>2$.  In particular, it can be generated by the following elements:

\begin{itemize}

\item The $n$ generators equivalent to the ones of $F_n$ described by \cite{Brown} (see \autoref{Theorem_Gen1}).

\item $5$ (for $n=2$) or $\frac{n^2}{2}+\frac{n}{2}+1$ (for $n> 2$) generators having the form $(T,A_{i,j},\lambda_{id},T)$, where~$A_{i,j}$ is any $n$-irreducible generator in some $\mathcal{P}_m$ satisfying $j-i=n$ or $m=n$, and $T$ can be any full $n$-ary tree with $m$ leaves. One can choose a different $T$ for each generator.

\item The $\dfrac{(n-1)n^2}{2}$  generators that can be written as $(R,1,\lambda_{k,A_{i,j}},R)$, for $0<i\leq n$, where $R$ is the full $n$-ary tree with one caret, and $\lambda_{k,A_{i,j}}=\{\ell_1,\dots,\ell_n\}$ is such that $\ell_k=A_{i,j}$ is a generator of $\mathcal{P}_n$ and $\ell_{k'}=1$ if $k'\neq k$.
\end{itemize}

\end{theorem}

\begin{figure}[h]
\centering
\includegraphics[scale=0.8]{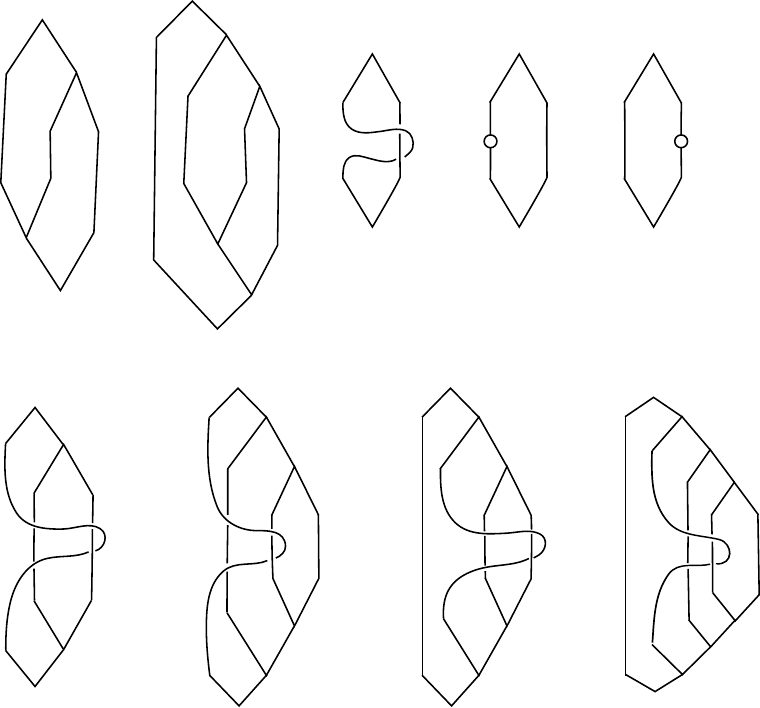}
\caption{Generators of $BF_2(\mathcal{P}_2)$. The label in the generators is $\sigma_1^{-2}$, the only generator of $\mathcal{P}_2$.}\label{Fig_gen2}
\end{figure}

\begin{question}
\emph{Which is the least number of elements that we can use to generate $BF_n(\{1\})$ and $BF_n(\mathcal{P}_n)$? }
\end{question}

\begin{question}
\emph{Can we find explicit presentations for $BF_n(\{1\})$ and $BF_n(\mathcal{P}_n)$? Can we find presentations with a minimal number of generators?}
\end{question}

\begin{conjecture}\label{conjecture_iso}
\emph{$BF_n(\{1\})$ and $BF_n(\mathcal{P}_n)$ are not isomorphic.}
\end{conjecture}

\begin{proof}[Proof of \autoref{Theorem_Gen3}]
For the number of the third item, just notice that the number of generators of $\mathcal{P}_n$ is~$\frac{(n-1)n}{2}$.

Let $R$ be the full $n$-ary tree with one caret.  Suppose that $A_{i,j}$ is an irreducible generator in~$P_{m}$ with $m>n$ and $j-i <n$, and consider $(T,A_{i,j},\lambda_{id},T)$. We claim that we can generate it using only generators from the statement. First observe that if $T'$ is any tree having $m-n+1$ leaves, then doing an expansion we have that 
$$  x:=(T',1,\lambda_{k,A_{i',j'}},T')\equiv
(T'[k],A_{i,j},\{1,\dots,1,\underbrace{A_{i',j'},\dots,A_{i',j'}}_{n \text{ (from position } i \text{ to } j)},1,\dots,1\},T'[k]),$$ for some $k,i'$ and $j'$. We claim that any element of this kind can be generated using only elements from the statement. We also claim the same for $$y:=(T'[k],1,\{1,\dots,1,\underbrace{A_{i',j'}^{-1},\dots,A_{i',j'}^{-1}}_{n \text{ (from position } i \text{ to } j)},1,\dots,1\},T'[k]).$$ If this is true, $(T,A_{i,j},\lambda_{id},T)$ is equivalent to the product $(T,1,\lambda_{id},T'[k])\cdot x\cdot y \cdot(T'[k],1,\lambda_{id},T)$.  

We need to show that we can generate any $(T'',1,\lambda_{k',A_{i,j}^{- 1}},T'')$ for $T''\in\{T',T'[k]\}$ using the generators of the statement. The proof is the same as in \autoref{Theorem_Gen2}: if we want $k'$ to be bigger than $n$ by expanding the first strand $q$ times, we can obtain from $(R,1,\lambda_{k,A_{i,j}^{\pm 1}},R)$, with $1<i\leq n$, a $(R',1,\lambda_{k+q(n-1),A_{i,j}^{\pm 1}},R')$ where $R'$ has $(q+1)n-q$ leaves. Otherwise, one can just expand the last leaf. We can repeat this process as many times as necessary and then multiply by $(T'',1,\lambda_{id}, R')$ and $(R',1,\lambda_{id}, T'')$ to obtain the desired element. 

To finish the proof, we need to count how many generators are in the second item of the proof. For $n=2$, this is just a finite checking. For $n>2$, we apply the same counting strategy as in the proof of \autoref{Theorem_Gen1}. Recall that, for $A_{i,j}$ to be irreducible $m\in\{n, 2n-1, 3n-2\}$. If $n=m$, then the number of generators is $\frac{(n-1)n}{2}$. For $m=2n-1$, we have $n-1$ possibilities, and for $m=3n-2$ we have two possibilities. In total this gives us $\frac{n^2}{2}+\frac{n}{2}+1$ elements.



\end{proof}

\begin{remark}
Notice that if we replace $H$ by $\mathcal{P}_n$ in \autoref{Theorem_Gen2}, so $k=\frac{(n-1)n}{2}$, we obtain $\frac{n^3}{2}+\frac{n^2}{2}+n+1$, so \autoref{Theorem_Gen3} reduces the number of generators by $\frac{(n-1)n}{2}$. 

Moreover, it is interesting to see that the smallest number of generator that we have found for $BF_2(\{1\})$ is 10, and for $BF_2(\mathcal{P}_2)$ is 9 (see those generators in \autoref{Fig_gen2}). This leads to some questions and results.
\end{remark}

\begin{question}
\emph{Which is the least number of elements that we can use to generate $BF_n(\{1\})$ and $BF_n(\mathcal{P}_n)$? }
\end{question}

\begin{question}
\emph{Can we find explicit presentations for $BF_n(\{1\})$ and $BF_n(\mathcal{P}_n)$? Can we find presentations with a minimal number of generators?}
\end{question}

\begin{proposition}\label{conjecture_iso}
$BF_n(\{1\})$ and $BF_n(\mathcal{P}_n)$ are not isomorphic.
\end{proposition}

\begin{proof}[Proof (based on ideas of Matthew Zaremsky).]
Notice that $BF_n(\mathcal{P}_n)$ surjects onto $BF_n \times \mathcal{P}_n$ by sending $(T_1, p, \{\lambda_1,\dots\}, T_2)$ to $(T_1,p, T_2) \cdot \lambda_1$. Let $\mathbb{Z}^k$ be the abelianization of $BF_n$ ---we know that the abelianization of $BF_2$ is $\mathbb{Z}^4$ \citep{Zaremsky2018}---. Also $\mathcal{P}_n$ surjects onto $\mathcal{P}_2\simeq \mathbb{Z}$ by forgetting strands. Then $BF_n(\mathcal{P}_n)$ surjects onto $\mathbb{Z}^{k+1}$ and therefore it cannot be isomorphic to~$BF_n$.
\end{proof}

\bigskip

\noindent{\textbf{\Large{Acknowledgments}}}

The author was supported by a Ram\'on y Cajal 2021 grant and the research grant PID2022-138719NA-I00 (Proyectos de Generaci\'on de Conocimiento 2022), both financed by the Spanish Ministry of Science and Innovation. She warmly thanks Matt Zaremsky for his interest in this work, providing the references for \autoref{remark_Ishida} and the proof of \autoref{conjecture_iso}.

\medskip

\medskip
\bibliography{Bibliography}

\bigskip\bigskip{\footnotesize%

\noindent
\textit{\textbf{María Cumplido} \\ 
Departmento de \'Algebra,
Facultad de Matem\'aticas,
Universidad de Sevilla. \\
Calle Tarfia s/n
41012, Seville, Spain.} \par
 \textit{E-mail address:} \texttt{\href{mailto:cumplido@us.es}{cumplido@us.es}}
 }

\end{document}

%% file: thompsonequiv.pdf_tex
\begingroup%
  \makeatletter%
  \providecommand\color[2][]{%
    \errmessage{(Inkscape) Color is used for the text in Inkscape, but the package 'color.sty' is not loaded}%
    \renewcommand\color[2][]{}%
  }%
  \providecommand\transparent[1]{%
    \errmessage{(Inkscape) Transparency is used (non-zero) for the text in Inkscape, but the package 'transparent.sty' is not loaded}%
    \renewcommand\transparent[1]{}%
  }%
  \providecommand\rotatebox[2]{#2}%
  \newcommand*\fsize{\dimexpr\f@size pt\relax}%
  \newcommand*\lineheight[1]{\fontsize{\fsize}{#1\fsize}\selectfont}%
  \ifx\svgwidth\undefined%
    \setlength{\unitlength}{388.6131563bp}%
    \ifx\svgscale\undefined%
      \relax%
    \else%
      \setlength{\unitlength}{\unitlength * \real{\svgscale}}%
    \fi%
  \else%
    \setlength{\unitlength}{\svgwidth}%
  \fi%
  \global\let\svgwidth\undefined%
  \global\let\svgscale\undefined%
  \makeatother%
  \begin{picture}(1,0.46868023)%
    \lineheight{1}%
    \setlength\tabcolsep{0pt}%
    \put(0,0){\includegraphics[width=\unitlength,page=1]{thompsonequiv.pdf}}%
    \put(0.05176993,0.39845382){\makebox(0,0)[lt]{\lineheight{1.25}\smash{\begin{tabular}[t]{l}$T_1$\end{tabular}}}}%
    \put(-0.00111324,0.20242548){\makebox(0,0)[lt]{\lineheight{1.25}\smash{\begin{tabular}[t]{l}$(1,2)(3,4)$\end{tabular}}}}%
    \put(0.09104453,0.02017498){\makebox(0,0)[lt]{\lineheight{1.25}\smash{\begin{tabular}[t]{l}$T_2$\end{tabular}}}}%
    \put(0,0){\includegraphics[width=\unitlength,page=2]{thompsonequiv.pdf}}%
    \put(0.92959602,0.39845382){\makebox(0,0)[lt]{\lineheight{1.25}\smash{\begin{tabular}[t]{l}$T_1[3]$\end{tabular}}}}%
    \put(0.9204625,0.20242548){\makebox(0,0)[lt]{\lineheight{1.25}\smash{\begin{tabular}[t]{l}$(1,2)(6,7)$\end{tabular}}}}%
    \put(0.88395391,0.02017498){\makebox(0,0)[lt]{\lineheight{1.25}\smash{\begin{tabular}[t]{l}$T_2[3]$\end{tabular}}}}%
    \put(0,0){\includegraphics[width=\unitlength,page=3]{thompsonequiv.pdf}}%
  \end{picture}%
\endgroup%

%% file: generador_puro.pdf_tex
\begingroup%
  \makeatletter%
  \providecommand\color[2][]{%
    \errmessage{(Inkscape) Color is used for the text in Inkscape, but the package 'color.sty' is not loaded}%
    \renewcommand\color[2][]{}%
  }%
  \providecommand\transparent[1]{%
    \errmessage{(Inkscape) Transparency is used (non-zero) for the text in Inkscape, but the package 'transparent.sty' is not loaded}%
    \renewcommand\transparent[1]{}%
  }%
  \providecommand\rotatebox[2]{#2}%
  \newcommand*\fsize{\dimexpr\f@size pt\relax}%
  \newcommand*\lineheight[1]{\fontsize{\fsize}{#1\fsize}\selectfont}%
  \ifx\svgwidth\undefined%
    \setlength{\unitlength}{233.39428999bp}%
    \ifx\svgscale\undefined%
      \relax%
    \else%
      \setlength{\unitlength}{\unitlength * \real{\svgscale}}%
    \fi%
  \else%
    \setlength{\unitlength}{\svgwidth}%
  \fi%
  \global\let\svgwidth\undefined%
  \global\let\svgscale\undefined%
  \makeatother%
  \begin{picture}(1,0.45280141)%
    \lineheight{1}%
    \setlength\tabcolsep{0pt}%
    \put(0,0){\includegraphics[width=\unitlength,page=1]{generador_puro.pdf}}%
    \put(0.25218009,0.41384963){\makebox(0,0)[lt]{\lineheight{1.25}\smash{\begin{tabular}[t]{l}$i$\end{tabular}}}}%
    \put(0.71266698,0.41766792){\makebox(0,0)[lt]{\lineheight{1.25}\smash{\begin{tabular}[t]{l}$j$\end{tabular}}}}%
  \end{picture}%
\endgroup%

%% file: generalizedelement.pdf_tex
\begingroup%
  \makeatletter%
  \providecommand\color[2][]{%
    \errmessage{(Inkscape) Color is used for the text in Inkscape, but the package 'color.sty' is not loaded}%
    \renewcommand\color[2][]{}%
  }%
  \providecommand\transparent[1]{%
    \errmessage{(Inkscape) Transparency is used (non-zero) for the text in Inkscape, but the package 'transparent.sty' is not loaded}%
    \renewcommand\transparent[1]{}%
  }%
  \providecommand\rotatebox[2]{#2}%
  \newcommand*\fsize{\dimexpr\f@size pt\relax}%
  \newcommand*\lineheight[1]{\fontsize{\fsize}{#1\fsize}\selectfont}%
  \ifx\svgwidth\undefined%
    \setlength{\unitlength}{189.85904003bp}%
    \ifx\svgscale\undefined%
      \relax%
    \else%
      \setlength{\unitlength}{\unitlength * \real{\svgscale}}%
    \fi%
  \else%
    \setlength{\unitlength}{\svgwidth}%
  \fi%
  \global\let\svgwidth\undefined%
  \global\let\svgscale\undefined%
  \makeatother%
  \begin{picture}(1,1.23046864)%
    \lineheight{1}%
    \setlength\tabcolsep{0pt}%
    \put(0,0){\includegraphics[width=\unitlength,page=1]{generalizedelement.pdf}}%
    \put(0.21435079,1.08546085){\makebox(0,0)[lt]{\lineheight{1.25}\smash{\begin{tabular}[t]{l}$T_1$\end{tabular}}}}%
    \put(0.29474004,0.04255981){\makebox(0,0)[lt]{\lineheight{1.25}\smash{\begin{tabular}[t]{l}$T_2$\end{tabular}}}}%
    \put(0,0){\includegraphics[width=\unitlength,page=2]{generalizedelement.pdf}}%
    \put(-0.00262325,0.67416041){\makebox(0,0)[lt]{\lineheight{1.25}\smash{\begin{tabular}[t]{l}$\sigma_4^2 \sigma_2\sigma_1^2\sigma_2$\end{tabular}}}}%
  \end{picture}%
\endgroup%

%% file: generalizedelement2.pdf_tex
\begingroup%
  \makeatletter%
  \providecommand\color[2][]{%
    \errmessage{(Inkscape) Color is used for the text in Inkscape, but the package 'color.sty' is not loaded}%
    \renewcommand\color[2][]{}%
  }%
  \providecommand\transparent[1]{%
    \errmessage{(Inkscape) Transparency is used (non-zero) for the text in Inkscape, but the package 'transparent.sty' is not loaded}%
    \renewcommand\transparent[1]{}%
  }%
  \providecommand\rotatebox[2]{#2}%
  \newcommand*\fsize{\dimexpr\f@size pt\relax}%
  \newcommand*\lineheight[1]{\fontsize{\fsize}{#1\fsize}\selectfont}%
  \ifx\svgwidth\undefined%
    \setlength{\unitlength}{206.27161942bp}%
    \ifx\svgscale\undefined%
      \relax%
    \else%
      \setlength{\unitlength}{\unitlength * \real{\svgscale}}%
    \fi%
  \else%
    \setlength{\unitlength}{\svgwidth}%
  \fi%
  \global\let\svgwidth\undefined%
  \global\let\svgscale\undefined%
  \makeatother%
  \begin{picture}(1,1.13256296)%
    \lineheight{1}%
    \setlength\tabcolsep{0pt}%
    \put(0,0){\includegraphics[width=\unitlength,page=1]{generalizedelement2.pdf}}%
    \put(0.61603203,1.04473098){\makebox(0,0)[lt]{\lineheight{1.25}\smash{\begin{tabular}[t]{l}$T_1[4]$\end{tabular}}}}%
    \put(0.56031069,0.05427601){\makebox(0,0)[lt]{\lineheight{1.25}\smash{\begin{tabular}[t]{l}$T_2[4]$\end{tabular}}}}%
    \put(0,0){\includegraphics[width=\unitlength,page=2]{generalizedelement2.pdf}}%
  \end{picture}%
\endgroup%

%% file: composition.pdf_tex
\begingroup%
  \makeatletter%
  \providecommand\color[2][]{%
    \errmessage{(Inkscape) Color is used for the text in Inkscape, but the package 'color.sty' is not loaded}%
    \renewcommand\color[2][]{}%
  }%
  \providecommand\transparent[1]{%
    \errmessage{(Inkscape) Transparency is used (non-zero) for the text in Inkscape, but the package 'transparent.sty' is not loaded}%
    \renewcommand\transparent[1]{}%
  }%
  \providecommand\rotatebox[2]{#2}%
  \newcommand*\fsize{\dimexpr\f@size pt\relax}%
  \newcommand*\lineheight[1]{\fontsize{\fsize}{#1\fsize}\selectfont}%
  \ifx\svgwidth\undefined%
    \setlength{\unitlength}{431.81196078bp}%
    \ifx\svgscale\undefined%
      \relax%
    \else%
      \setlength{\unitlength}{\unitlength * \real{\svgscale}}%
    \fi%
  \else%
    \setlength{\unitlength}{\svgwidth}%
  \fi%
  \global\let\svgwidth\undefined%
  \global\let\svgscale\undefined%
  \makeatother%
  \begin{picture}(1,0.88042067)%
    \lineheight{1}%
    \setlength\tabcolsep{0pt}%
    \put(0,0){\includegraphics[width=\unitlength,page=1]{composition.pdf}}%
    \put(-0.00081352,0.82524753){\makebox(0,0)[lt]{\lineheight{1.25}\smash{\begin{tabular}[t]{l}$T_1$\end{tabular}}}}%
    \put(0.0278874,0.45290616){\makebox(0,0)[lt]{\lineheight{1.25}\smash{\begin{tabular}[t]{l}$T_2$\end{tabular}}}}%
    \put(0.58001289,0.86067181){\makebox(0,0)[lt]{\lineheight{1.25}\smash{\begin{tabular}[t]{l}$T_1[6]$\end{tabular}}}}%
    \put(0.49840502,0.43932899){\makebox(0,0)[lt]{\lineheight{1.25}\smash{\begin{tabular}[t]{l}$T_2[6]=T_3[2]$\end{tabular}}}}%
    \put(0.02885085,0.41232925){\makebox(0,0)[lt]{\lineheight{1.25}\smash{\begin{tabular}[t]{l}$T_3$\end{tabular}}}}%
    \put(0.02889247,0.03302685){\makebox(0,0)[lt]{\lineheight{1.25}\smash{\begin{tabular}[t]{l}$T_4$\end{tabular}}}}%
    \put(0.57834601,0.01303835){\makebox(0,0)[lt]{\lineheight{1.25}\smash{\begin{tabular}[t]{l}$T_4[2]$\end{tabular}}}}%
    \put(0.95135899,0.66257682){\makebox(0,0)[lt]{\lineheight{1.25}\smash{\begin{tabular}[t]{l}$T_1[6]$\end{tabular}}}}%
    \put(0.9508854,0.2302268){\makebox(0,0)[lt]{\lineheight{1.25}\smash{\begin{tabular}[t]{l}$T_4[2]$\end{tabular}}}}%
  \end{picture}%
\endgroup%